\documentclass[12pt]{article}
\usepackage{amssymb, amsmath, amsthm, verbatim, cite, lscape}
\usepackage[cp1251]{inputenc}
\usepackage[english]{babel}
\usepackage{srcltx,mathrsfs,textcomp}

\newtheorem{theorem}{Theorem}

\newcounter{parag}
\newcommand{\sect}[1]
{\refstepcounter{parag}
\begin{center} { \bf\theparag. #1} \end{center}}
\newtheorem{cor}{Corollary}
\newtheorem{lemma}{Lemma}[parag]
\newtheorem{prop}[lemma]{Proposition}
\theoremstyle{definition}
\newtheorem*{prf}{Proof}
\newtheorem*{theorem*}{Theorem}

\theoremstyle{definition}
\newtheorem*{rem}{Remark}

\begin{document}

\begin{center} {\large Groups with bounded centralizer chains and the~Borovik--Khukhro conjecture}
\end{center}

\begin{center} A.\,A.\,Buturlakin\footnote[1]{Sobolev Institute of Mathematics, 4 Acad. Koptyug avenue, Novosibirsk 630090, Russia}$^,$\footnote[2]{Novosibirsk State University, 1, Pirogova str., Novosibirsk, 630090, Russia}$^,$\footnote[3]{Corresponding author

\emph{E-mail addresses:} buturlakin@math.nsc.ru, revin@math.nsc.ru, vasand@math.nsc.ru}, D.\,O.\,Revin\footnotemark[1]$^,$\footnotemark[2], A.\,V.\,Vasil$'$ev\footnotemark[1]$^,$\footnotemark[2]

\end{center}

{\small\begin{center}\textbf{Abstract} \end{center}

Let $G$ be a locally finite group and $F(G)$ the Hirsch--Plotkin radical of $G$. Denote by $S$ the full inverse image of the generalized Fitting subgroup of $G/F(G)$ in $G$. Assume that there is a number $k$ such that the length of every chain of nested centralizers in $G$ does not exceed $k$.  The Borovik--Khukhro conjecture states, in particular, that under this assumption the quotient $G/S$ contains an abelian subgroup of index bounded in terms of $k$. We disprove this statement and prove some its weaker analog.

}

\begin{center}\textbf{Introduction}\end{center}

Following \cite{MS}, the maximal length of chains of nested centralizers of a group $G$ is called the \emph{$c$-dimension} of $G$ (the same number is also known as the length of the centralizer lattice of $G$, see, e.g.,~\cite{Schmidt}). Here we consider groups of finite $c$-dimension. The class of such groups is quite wide: it includes, for example, abelian groups, linear groups, torsion-free hyperbolic groups, and hence free groups. It is closed under taking subgroups, finite direct products, and extensions by finite groups; however, the $c$-dimension of a homomorphic image of a group from this class can be infinite.
Obviously, each group of finite $c$-dimension satisfies the minimal condition on centralizers, therefore, it is a so-called \emph{$\mathfrak M_C$-group}.

The minimal condition on centralizers is a useful notion in the class of locally finite groups. Indeed, locally finite $\mathfrak M_C$-groups satisfy an analog of the Sylow theorem \cite[Teorem~B]{Bryant}. Periodic locally soluble $\mathfrak M_C$\nobreakdash-groups satisfy an analog of the Hall theorem on $\pi$-subgroups \cite[Theorem~1.6]{BH}, moreover, these groups are soluble and nilpotent-by-abelian-by-finite \cite[Theorems~2.1 and 2.2]{BH}. When we restrict ourselves to the class of locally finite groups of finite $c$-dimension, it is natural to ask if some strengthening of mentioned results can be proved. For example, whether the derived length of a locally soluble groups from this class is bounded in terms of $c$-dimension? Another natural question is whether such a group contains a nilpotent-by-abelian subgroup of index bounded in the same terms. The first question was answered in the affirmative \cite[Theorem~1]{Khukhro}. The answer to the second question is ``no'' \cite[Example~1]{Khukhro}, but a weaker statement holds: if $G$ is a periodic locally soluble group of $c$-dimension $k$, then the quotient group of $G$ by the second Hirsch--Plotkin radical $F_2(G)$ contains an abelian subgroup of $k$-bounded index.

In the same paper \cite{Khukhro} E.I.\,Khukhro posed the following conjecture which he attributed to A.V.\,Borovik.

\smallskip

\textbf{The Borovik--Khukhro conjecture}. {\it Let $G$ be a locally finite group of finite
$c$-dimension $k$. Let $S$ be the full inverse image of $F^*(G/F(G))$ in $G$. Then

$(1)$ the number of nonabelian simple composition factors of $G$ is finite and
$k$-bounded;

$(2)$ $G/S$ has an abelian subgroup of finite $k$-bounded index.}

\smallskip

Recall that the Hirsch--Plotkin radical $F(G)$ is the maximal normal locally nilpotent subgroup. The $i$-th Hirsch--Plotkin radical $F_i(G)$ is defined inductively: $F_1(G)=F(G)$ and $F_i(G)$ is the full inverse image of $F(G/F_{i-1}(G))$ in~$G$. The layer $E(G)$ of $G$ is the subgroup generated by all components of~$G$, that is its subnormal quasisimple subgroups. The generalized Fitting subgroup $F^*(G)$ is the product of $F(G)$ and $E(G)$.

While the part $(1)$ of the conjecture is known to be true \cite{13ButVas}, the second one turns out to be false.

\begin{theorem}\label{t:ConjIsFalse} The statement $(2)$ of the Borovik--Khukhro conjecture does not hold.
\end{theorem}

It is clear that a counterexample proving Theorem \ref{t:ConjIsFalse} should be an infinite series of locally finite groups $G$ of uniformly bounded $c$-dimensions with unbounded minimal indices of abelian subgroups in the quotients $G/S$. In fact, we construct an infinite number of counterexamples as follows. For each positive integer $n$, we define an infinite series $\mathscr{G}_n$ of finite groups $G_{n,r}$ where $r$ runs over an infinite set of primes. The $c$-dimensions of groups in $\mathscr{G}_n$ are bounded in terms of $n$ only. The quotient of $G=G_{n,r}$ by the corresponding subgroup $S$ (that is by the full inverse image of generalized Fitting subgroup of $G_{n,r}/F(G_{n,r})$) is isomorphic to the symplectic group $Sp_{2n}(r)$. This contradicts the statement $(2)$ of the Borovik--Khukhro conjecture, because the minimum of indices of proper subgroups of $G/S=Sp_{2n}(r)$ grows as $r$ increases.

It is worth noting that our construction depends on some deep number-theoretical results. It requires some results on Dickson's conjecture and Schinzel's hypothesis H on prime values of finite sets of linear forms and irreducible integer-valued polynomials obtained in \cite{AST}. We are thankful to participants of the forum on mathoverflow.net, especially to T. Tao, and to J. Maynard for the information on various results on these conjectures and related seive methods. There is more detailed discussion on this topic after the proof of Theorem \ref{t:ConjIsFalse} in Section 2.

Though the second part of the Borovik--Khukhro conjecture is false, some its weaker analog holds at least in the case of finite groups.

\begin{theorem}\label{t:3Fitting} Let $G$ be a finite group of finite $c$-dimension $k$. Put $\overline G=G/F_3(G)$. The group $\overline G/E(\overline G)$ contains an abelian subgroup of $k$-bounded index.
\end{theorem}

The proof of Theorem \ref{t:3Fitting} relies on \cite[Theorem 1]{Khukhro} and ideas from \cite{13ButVas}. The difference from the proof in \cite{13ButVas} is that we bound not just the number of nonabelian composition factors, but the common ``size'' of them.


Theorem \ref{t:3Fitting} implies that, for a finite group $G$, the quotient of $G/F_3(G)$ by the generalized Fitting subgroup contains an abelian subgroup of index bounded in term of $c$-dimension, while it was conjectured that $F(G)$ instead of $F_3(G)$ is enough. So the natural question is whether one can replace $F_3(G)$ with $F_2(G)$. Furthermore, Theorem \ref{t:3Fitting} regards only finite groups. Hence it would be interesting to obtain some positive statement for locally finite groups. Since every iterated Hirsch--Plotkin radical is contained in a locally soluble radical, the following seems to be the ``weakest'' analog of the conjecture: Is it true that the quotient group of $G/R$ by the generalized Fitting subgroup contains an abelian subgroup whose index is bounded in terms of $k$ (here $G$ is locally finite of $c$-dimension $k$ and $R$ is a (locally) soluble radical of $G$)? The affirmative answer to this question for finite groups follows from Theorem \ref{t:3Fitting} (see Corollary \ref{t:solubleRadical} in Section 3).

Section 1 contains preliminary results on group centralizers and some information on simple groups. In Section 2, we give a counterexample to the Borovik--Khukhro conjecture and discuss for which $c$-dimensions such counterexamples exist. Theorem \ref{t:3Fitting} is proved in Section 3.

We thank E.I. Khukhro for drawing our attention to the subject.

\sect{Preliminaries}

The following lemma lists some well-known properties of centralizers.

\begin{lemma}\label{l:PropertyOfCentralizers} Let $G$ be a group and $X$, $Y$ be subsets in $G$. Then
\begin{itemize}
 \item[$(1)$] $C_G(X)=\bigcap\limits_{x\in X}C_G(x)$;
 \item[$(2)$] $C_G(X)=C_G(\langle X\rangle)$;
 \item[$(3)$]$C_G(X)\leqslant C_G(Y)$  if and only if  $C_G(C_G(X))\geqslant C_G(C_G(Y))$; furthermore, $C_G(X)=C_G(Y)$  if and only if  $C_G(C_G(X))=C_G(C_G(Y))$.
\end{itemize}
\end{lemma}



Let $G$ be a group and $V$ a normal abelian subgroup of $G$. Let $\overline{\phantom{g}}:G\rightarrow G/V$ be the natural homomorphism. The action of $G$ on $V$ by conjugation induces the action of $G/V$ on $V$. In the following lemma the centralizer $C_V(\overline Y)$ for $Y\leqslant G$ is defined with respect to this action, and therefore coincides with~$C_V(Y)$.



\begin{lemma}\label{l:c-dimensionOfSemidirectProduct} Let $G$ be a group and $V$ a normal abelian subgroup of $G$. Suppose that
 the length of every chain of nested subgroups in $G/V$ is at most~$l$. Then the $c$-dimension of $G$ is at most~$2l$. If $V$ is central, then $c$-dimension of $G$ is at most~$l$.
\end{lemma}

\begin{prf} Suppose that the $c$-dimension of $G$ is equal to $k$. Let\footnote{Here and further, $H<K$ denotes the strict inclusion, that is $H\leqslant K$ and $H\neq K$.}
\begin{equation} \label{e0}
 Y_0>\dots> Y_k
\end{equation}
be a chain of subgroups of $G$ such that
\begin{equation} \label{e1}
C_G(Y_0)<\dots< C_G(Y_k).
\end{equation}

As above $\overline{\phantom{g}}:G\rightarrow G/V$ is the natural homomorphism. In view of (\ref{e0}) and (\ref{e1}), we have the following chains of inclusions:

\begin{equation} \label{e1.5}
\overline Y_0 \geqslant\dots\geqslant \overline Y_k ,
\end{equation}
\begin{equation} \label{e2}
C_V(Y_0)\leqslant\dots\leqslant C_V(Y_k),
\end{equation}
\begin{equation} \label{e3}
\overline {C_G(Y_0)} \leqslant\dots\leqslant \overline {C_G(Y_k)},
\end{equation}
\begin{equation} \label{e4}
C_V(\overline Y_0)\leqslant\dots\leqslant C_V(\overline Y_k).
\end{equation}
Since
$C_G(Y_i)/C_V(Y_i)\cong \overline{C_G(Y_i)}$  and $C_G(Y_{i-1})< G_G(Y_{i})$ for every $i$, at least one
of the inclusions  $$C_V(\overline Y_{i-1})=C_V(Y_{i-1})\leqslant C_V(Y_{i})= C_V(\overline Y_{i}) \text{ and
} \overline{C_G(Y_{i-1})} \leqslant \overline{C_G(Y_{i})}$$ is strict. Moreover, if $C_V(\overline Y_{i-1})< C_V(\overline Y_{i})$, then $\overline Y_{i-1} > \overline Y_{i}$.
Thus, $k$ does not exceed the sum of the numbers of strict inclusions in the chains~(\ref{e1.5}) and~(\ref{e3}) of subgroups of~$G/V$. Hence $k\leqslant 2l$. Finally, if $V$ is cental then all inclusions in $(4)$ are equalities, and all inclusions in $(5)$ are strict.
\end{prf}

\begin{rem} We want to point out that the bound for $c$-dimension obtained in Lemma \ref{l:c-dimensionOfSemidirectProduct} for the general case is sharp. Let $H$ be a finite group and $$1=H_0<H_1<\dots<H_l=H$$ be its longest chain of nested subgroups. Fix some prime $p$. Denote by $V_i$ for $i=0,\dots, l$ the permutation module of $H$ over the field of order $p$ given by the action of $H$ on the right cosets by $H_i$ via right multiplications. Put $$V=V_0\oplus V_1\oplus\dots\oplus V_l.$$ Denote by $G$ the natural semidirect product $VH$. We calim that the $c$-dimension of $G$ is equal to $2l$.

Let $e_i$ be a vector of $V_i$ such that $C_H(e_i)=H_i$. Let $W_i$ be the subspace of $V$ spanned by $e_i,e_{i+1},\dots, e_l$.  It is easy to see that $C_G(W_i)=VH_i$. Also the centralizer of $VH_i$ contains $W_i$ and does not contain $e_{i-1}$. Therefore we have $$C_G(W_l)>C_G(W_{l-1})>\dots>C_G(W_0)=C_G(V)$$$$=C_G(VH_0)>C_G(VH_1)>\dots>C_G(VH_l)=Z(G).$$ The length of this chain is exactly $2l$, i.e. maximum possible in terms of the lemma.
\end{rem}


The following lemma is crucial to the proof of Theorem \ref{t:3Fitting}.

\begin{lemma}\label{LemmaKhukhro}{\rm \cite[Lemma 3]{Khukhro}} If an elementary abelian $p$-group $E$ of order~$p^n$ acts faithfully on a finite nilpotent $p'$-group $Q$, then there exists a series of subgroups $$E=E_0>E_1>E_2>\dots>E_n=1$$ such that $$C_Q(E_0)<C_Q(E_1)<\dots<C_Q(E_n).$$
\end{lemma}

\begin{lemma}\label{l:CentralProduct} The $c$-dimension of a central product of finite groups $H$ and $K$ is at least the sum of their $c$-dimensions.
\end{lemma}

\begin{prf} Let $G$ be a central product $H\circ K$. Let $$H_0<H_1<\dots<H_m \quad\text{ and }\quad K_0<K_1<\dots<K_l$$ be series of subgroups of $H$ and $K$ such that the series of their centralizers are strict and have maximal lengths. We may assume that both $H$ and $K$  are subgroups of $G$, and, therefore, all their subgroups are also subgroups of $G$. We have $C_G(H_i)=C_H(H_i)\circ K$ and $C_G(HK_i)=Z(H)\circ C_K(K_i)$. Hence the subgroups $H_i$ and $HK_i$ provide a series of subgroups with centralizer chain of desired length.
\end{prf}

\begin{rem} The $c$-dimension of $H\circ K$ can be larger than the sum of $c$\nobreakdash-dimensions of $H$ and $K$. Indeed, let $H$ be the group $SL_2(3)$. Then $H$ is a semidirect product of the quaternion group of order $8$ and the cyclic group of order $3$ which permutes subgroups of order $4$. It is not hard to see that the $c$-dimension of $H$ is $2$.  Denote by $Q$ the Sylow $2$-subgroup of $H$. If $Q$ is generated by two elements $a$ and $b$, then the commutator $[a,b]$ is the generator of the center of $Q$ and, therefore, of $H$. Let $K$ be a group isomorphic to $H$ and $c$ and $d$ be the images of $a$ and $b$ under this isomorphism.

Put $G=(H\times K)/\langle [a,b][c,d]\rangle$ and denote by $\overline{\phantom{g}}$ the canonical homomorphism from $H\times K$ onto $G$. By definition, $G$ is a central product of $H$ and~$K$. Consider the following series of centralizers:\begin{equation}\label{e:Remark}G>C_G(\bar a)\geqslant C_G(\bar a\bar c)\geqslant C_G(\bar a,\bar c)>C_G(\bar H,\bar c)>Z(G).\end{equation}

We have $\bar b\bar d\in C_G(\bar a\bar c)\setminus C_G(\bar a,\bar c)$. Indeed, $\bar a^{\bar b\bar d}=\bar a^{-1}$ and $[\bar a\bar c,\bar b\bar d]=[\bar a,\bar b][\bar c,\bar d]=1$. Therefore the third inclusion in (\ref{e:Remark}) is strict.

Let us show that no element of order $3$ centralizes $\bar a\bar c$. Let $\bar g$ be an element of $G$ of order $3$. Without loss of generality, we may assume that $\bar g$ does not normalize any subgroup of order $4$ in $H$. Moreover, we may assume that $\bar a^{\bar g}=\bar b$. Suppose that $(\bar a\bar c)^{\bar g}=\bar a\bar c$. Hence $(\bar b\bar a^{-1})(\bar c^{\bar g}\bar c^{-1})=1$, but $ba^{-1}$ is an element of order $4$ in $H$ and its product with any $2$-element $k$ of $K$ is again an element of order $4$. Therefore $\bar b\bar a^{-1}\bar k$ is never equal to $1$. Hence $\bar g$ does not centralize $\bar a\bar c$.

Thus all inclusions in (\ref{e:Remark}) are strict and the $c$-dimension of $H\circ K$ is at least $5$ (direct calculations show that it is exactly $5$).
\end{rem}

For a finite group $G$, we define a non-negative integer $\lambda(G)$ as follows. First, let $G$ be a nonabelian simple group. If $G$ is a group of Lie type, then $\lambda(G)$ is the minimum of Lie ranks of groups of Lie type isomorphic to $G$ (the Lie rank is the rank of the corresponding $(B,N)$-pair \cite[page 249]{Asch}). If $G$ is an alternating group, then let $\lambda(G)$ be the degree of $G$ (except for the groups $A_5$, $A_6$ and $A_8$ which are isomorphic to groups of Lie type and, therefore, already have assigned values of $\lambda(G)$). Put $\lambda(G)=1$ for the sporadic groups. Now if $G$ is an arbitrary finite group, then $\lambda(G)$ is the sum of $\lambda(S)$ where $S$ runs over the nonabelian composition factors of $G$.

Recall that the $p$-rank of a group $G$ for a given prime $p$ is the largest rank of elementary abelian $p$-subgroups of $G$.

\begin{lemma}\label{l:ElemAbelian} Let $r\in\{2,3,5\}$. There exist a constant $c>0$ such that the $r$-rank of any nonabelian finite simple group $G$ whose order is divisible by $r$ is at least $c\lambda(G)$.
\end{lemma}

\begin{prf} For the  sporadic groups, there is nothing to prove. In the case of alternating groups, this fact is obvious. In the case of groups of Lie type, it follows, for example, from \cite[Part I, 10-6]{83GorLyo}.
\end{prf}

\begin{lemma}\label{l:NumberOfPrimeDivisors} Let $G$ be a finite simple group of Lie type. The number of distinct prime divisors of the order of $G$ is at least $\lambda(G)$.
\end{lemma}

\begin{prf} Lemma follows from formulas for orders of groups of Lie type \cite[Table~16.1]{Asch} and the well-known Zsigmondy theorem \cite{Zs}.
\end{prf}









Denote by $\mu(G)$ the degree of the minimal faithful permutation representation of a finite group~$G$.

\begin{lemma}\label{DegBySolRad}{\rm \cite[Theorem 2]{Holt}} Let $G$ be a finite group. Let $\mathfrak{L}$ be a class of finite groups closed under taking subgroups, homomorphic images and extensions. Let $N$ be the $\mathfrak{L}$-radical of $G$, that is the maximal normal $\mathfrak{L}$-subgroup of $G$. Then $\mu(G)\geqslant\mu(G/N)$.
\end{lemma}

\begin{lemma}\label{DirProdSimple}{\rm \cite[Theorem 3.1]{Praeger}} Let $S_1$, $S_2$, $\dots$, $S_r$ be simple groups. Then $\mu(S_1\times S_2\times\dots\times S_r)=\mu(S_1)+\mu(S_2)+\dots+\mu(S_r)$.
\end{lemma}

\begin{lemma}\label{l:MinPermPres} If $S$ is a finite nonabelian simple group, then $\mu(S)\geqslant \lambda(S)$.
\end{lemma}

\begin{prf} For alternating and sporadic groups, the assertion is obvious. If $S$ is a group of Lie type having a faithful permutation representation of degree $n$, then the number of prime divisors of $|S|$ is less than $n$ and, by Lemma \ref{l:NumberOfPrimeDivisors}, is not less than $\lambda(S)$.
\end{prf}

Note that Lemma \ref{l:MinPermPres} can also be deduced from the information on the values of $\mu(S)$ for finite simple groups $S$ of Lie type which are all known (the complete list of these numbers can be found, for example, in \cite[Table 4]{GMPS}).

\begin{lemma}\label{CompFactorSn} If $G$ is a subgroup of the symmetric group $\operatorname{Sym}_n$, then ${\lambda(G)<5n/4}$.
\end{lemma}

\begin{prf} Proceed by induction on $|G|$. If the soluble radical $R$ of $G$ is non-trivial, then Lemma \ref{DegBySolRad} implies that $G/R$ is also a subgroup of $\operatorname{Sym}_n$ and the inequality follows by induction.

Let $R$ be trivial. If the socle $\operatorname{Soc}(G)$ of $G$ is the direct product of nonabelian simple groups $S_1$, $S_2$, $\dots$, $S_l$, then $G$ is a subgroup of the semidirect product of ${\operatorname{Aut}}(S_1)\times {\operatorname{Aut}}(S_2)\times\dots\times {\operatorname{Aut}}(S_l)$ and some subgroup of $\operatorname{Sym}_l$. Lemmas~\ref{DirProdSimple} and  \ref{l:MinPermPres} imply that $\lambda(\operatorname{Soc}(G))\leqslant n$. Since $\mu(S)\geqslant 5$ for every finite simple group $S$, it follows from Lemma~\ref{DirProdSimple} that $l\leqslant n/5$. By the Schreier conjecture, the natural homomorphism from $G/\operatorname{Soc}(G)$ to $\operatorname{Sym}_l$ has soluble kernel. By inductive hypothesis $\lambda(G/\operatorname{Soc}(G))\leqslant 5l/4\leqslant n/4$. Finally, we have  $$\lambda(G)=\lambda(\operatorname{Soc}(G))+\lambda(G/\operatorname{Soc}(G))\leqslant n+n/4=5n/4.$$
\end{prf}

\begin{lemma}\label{l:NumberOfCompFactors}{\rm \cite[Proposition 2.1]{13ButVas}} If $G$ is a finite group of $c$-dimension $k$, then the number of nonabelian composition factors of $G$ is less than $5k$.
\end{lemma}

The following lemma is a direct consequence of \cite[Theorem 1]{Khukhro}.

\begin{lemma}\label{l:KhukhroTh} If a periodic locally soluble group $G$ has finite $c$-dimension $k$, then the quotient $G/F_2(G)$ has an abelian subgroup of finite $k$-bounded index. In particular, the order of $G/F_3(G)$ is $k$-bounded.
\end{lemma}

\sect{Counterexamlpes to the Borovik--Khukhro conjecture}

Denote by $\Omega(n)$ the number of prime divisors of a positive integer $n$ counting multiplicities, i.e. if $n=p_1^{\alpha_1}\cdots p_m^{\alpha_m}$ for  primes $p_1,\dots,p_m$, then $$\Omega(n)={\alpha_1}+\cdots +{\alpha_m}.$$

For positive integers $n$ and $M$, we set
$$\pi_{n,M}=\{r\mid r \text{ is an odd prime and } \Omega(r^n-1)\leqslant M\}.$$
It is clear that the following inclusions hold for every $n$:
\begin{equation} \label{e5}
\pi_{n,1}\subseteq\pi_{n,2}\subseteq\pi_{n,3}\subseteq\dots.
\end{equation}
A crucial ingredient in constructing our counterexamples is the following number-theoretical statement.

\begin{prop}\label{p:NumberTheory} {\rm \cite[Theorem C]{AST}} For every positive integer $n$, there exists a positive integer $M$ such that $\pi_{n,M}$ is infinite.
\end{prop}

Given $n$,  put $$M_n=\min\{M\mid  \pi_{n,M}\text{ is infinite}\}\,\,\,\,\text{ and }\,\,\,\, \pi_n=\pi_{n,M_n}.$$ Thus $\pi_n$ is the first infinite set in the
chain~(\ref{e5}).

 Fix a positive integer $n$ and an odd prime $r$. Denote by  $$R_{n,r} \text{ an extra special group of order } r^{2n+1}  \text{ and of exponent } r.$$
It is known that
   $\operatorname{Aut}(R_{n,r})$ is split over $\operatorname{Inn}(R_{n,r})$ \cite[p.~404]{Griess} and
    the image in $\operatorname{Out}(R_{n,r})$ of the centralizer in  $\operatorname{Aut}(R_{n,r})$ of $Z(R_{n,r})$ is isomorphic to the symplectic group $Sp_{2n}(r)$ \cite[ex.~8.5, p.~116]{Asch}.
Hence $\operatorname{Aut}(R_{n,r})$ contains a subgroup $$A_{n,r}\cong Sp_{2n}(r)$$ which centralizes $Z(R_{n,r})$. We can form the natural semidirect product
$$X_{n,r}=R_{n,r}A_{n,r}.$$

Take a prime $p$ with $p\equiv 1\pmod r$. It is known \cite[p.~151]{KL} that $R_{n,r}$ has a faithful irreducible representation of degree $r^n$
over the field
$\mathbb{F}_p$ of order~$p$. Moreover, this representation extends to a faithful representation of~$X_{n,r}$ \cite[(3A) and (3B)]{AlF}.
Denote by
$$V_{n,r}\text{  a faithful irreducible  } \mathbb{F}_pX_{n,r}\text{-module of dimension } r^n$$  corresponding to this representaion
and  form the natural semidirect product
$$G_{n,r}=V_{n,r}X_{n,r}.$$

The following statement lists some basic properties of $G_{n,r}$ which can be easily deduced from the definition of this group.

\begin{prop}\label{p:ConstructionOfCounterexample}
\begin{itemize}
 \item[$(1)$] $F(G_{n,r})=V_{n,r}$;
 \item[$(2)$] $F^*(G_{n,r}/V_{n,r})=F(G_{n,r}/V_{n,r})=R_{n,r}$;
 \item[$(3)$] $G_{n,r}/R_{n,r}\cong A_{n,r}\cong Sp_{2n}(r)$.
\end{itemize}
\end{prop}

For a positive integer $n$, put $$\overline n= 2\operatorname{lcm}(1,2,\dots,n) \text{ and }
\mathscr{G}_n=\left\{G_{n,r}\mid r\in\pi_{\overline{n}}\right\}.$$ Observe that by definition the set $\mathscr{G}_n$ is infinite.

\begin{theorem}\label{t:counterexample1} For an arbitrary positive integer $n$, the $c$-dimension of every group in $\mathscr{G}_n$ does not exceed
\begin{equation} \label{e6}
   2\left((n+1)^2+nM_{\overline{n}}\right).
\end{equation}

\end{theorem}

\begin{prf}
Take $G=G_{n,r}\in\mathscr{G}_n$. By definition we have $r\in\pi_{\overline{n}}$ and hence $$\Omega(r^{\overline{n}}-1)\leqslant M_{\overline{n}}.$$
The Lagrange theorem implies that the length of every chain of  nested subgroups in $X_{n,r}$  does not exceed $l=\Omega(|X_{n,r}|)$. Since
$$
|X_{n,r}|=|R_{n,r}||A_{n,r}|=r^{2n+1}|Sp_{2n}(r)|=r^{(n+1)^2}\prod\limits_{i=1}^n(r^{2i}-1)
$$
and $r^{2i}-1$ divides $r^{\overline{n}}-1$ for every $i=1,\dots,n$, we have
$$l= (n+1)^2+\sum\limits_{i=1}^n\Omega(r^{2i}-1)\leqslant (n+1)^2+n\Omega(r^{\overline{n}}-1)\leqslant (n+1)^2+nM_{\overline{n}}.$$
Now the theorem follows from Lemma~\ref{l:c-dimensionOfSemidirectProduct}.
\end{prf}

Let us prove Theorem~\ref{t:ConjIsFalse}. Assume that the statement $(2)$ of the Borovik--Khukhro conjecture is true. Fix a positive integer $n$ and consider the set $\mathscr{G}_n$. Theorem  \ref{t:counterexample1} implies that the $c$-dimensions of groups in $\mathscr{G}_n$ are bounded by the number given in (\ref{e6}). Hence every group $G_{n,r}/R_{n,r}\cong Sp_{2n}(r)$ for $r\in\pi_{\overline n}$ contains a (proper) abelian subgroup of index bounded in terms of $n$, and, therefore, a normal abelian subgroup of bounded index. But there is a unique maximal proper normal subgroup of $Sp_{2n}(r)$ which is its center of order $2$ whose index grows as $r$ increases. This contradiction completes the proof of Theorem~\ref{t:ConjIsFalse}.

\smallskip


Theorem \ref{t:ConjIsFalse} implies that there exists a positive integer $k$ such that the second part of the Borovik--Khukhro conjecture fails for the finite groups of $c$-dimension $k$. On the other hand, the description of finite groups of $c$-dimension $2$ \cite[Theorem 9.3.12]{Schmidt} yields that the full version of the conjecture holds for these groups. So it is natural to ask for which minimal value of the $c$-dimension the Borovik--Khukhro conjecture is false. We already know that $3$ is a lower bound for this number and we do not know any better lower bound. Let us obtain an upper bound by bounding $M_2$.

Recall that $M_2$ is the smallest $M$ such that there exist an infinite number of prime numbers
$r$ satisfying $\Omega(r^2-1)\leqslant M.$

Since $r^2-1$ is divisible by $24$ for $r>3$, it is clear that $\Omega(r^2-1)\geqslant 6$ for sufficiently large $r$, and therefore $M_2\geqslant 6$. There is a question attributed to P. Neumann \cite[p.~316]{Solomon} asking whether the number of prime numbers $r$ such that the order of $PSL_2(r)$ is a product of six prime numbers is infinite. Since
$$\Omega(r^2-1)=\Omega(|PSL_2(r)|),$$
the positive answer to this question would imply that $M_2=6$ and the $c$-dimensions of groups of $\mathscr{G}_1$ would not exceed $2\cdot(2^2+6)=20$.

The affirmative answer to Neumann's question would follow from the validity of Dickson's conjecture (the $k$-tuple conjecture) \cite{Dickson}. This conjecture states that for a finite set of integer linear forms $a_1n+b_1,\dots, a_kn+b_k$ with $a_i>0$, there are infinitely many positive integers $n$ such that all these forms are simultaneously prime unless there are fixed divisors of their product. Equivalently $$\Omega\left(\prod\limits_{i=1}^k(a_in+b_i)\right)=k$$ for infinitely many integer values of $n$. If Dickson's conjecture is valid for the triplet  $12n+1$, $6n+1$ and $n$, then there are infinitely many primes $r$ of the form $r=12n+1$ such that  $n$ and $6n+1$ are also primes and
$$\Omega(r^2-1)=\Omega(12n\cdot (12n+2))=\Omega(2^3\cdot 3\cdot n\cdot (6n+1))=6.$$

Both Dickson's and Neumann's conjectures provide an exact value of $M_2$, but for our purposes a partial result might be satisfactory. There are well developed sieve methods that allow to obtain partial results for this kind of statements (see \cite{Sieve} for details). For example, the main result of \cite{Maynard} is very close to what is required. In this paper, J.\,Maynard proved that for any triplet of the forms $a_in+b_i$, $i=1,2,3$, without fixed divisors of their product, there are infinitely many positive integer values of $n$ such that $$\Omega\left(\prod\limits_{i=1}^3(a_in+b_i)\right)\leq 7.$$ In particular, there are an infinite number of integer values of $n$ such that the product $n\cdot(6n+1)\cdot(12n+1)$ has at most seven prime divisors, and, therefore, there are infinitely many integers $m$ such that $\Omega(m(m^2-1))\leq 11$. Unfortunately, one cannot guarantee the primality of $m$ in this result.  Nevertheless in a private letter, J. Maynard expressed confidence that one can show that there are an infinite number of primes $r$ such that $\Omega(r(r^2-1))\leq 11$, or equivalently  $\Omega(r^2-1)\leq 10$, by slightly changing his proof. This would mean that the $c$-dimensions of groups in $\mathscr{G}_1$ do not exceed $28$.

Proposition~\ref{p:NumberTheory}, which is also proved in \cite{AST} by using sieve methods, has the following refinement \cite[Corollary~4.2]{AST}: there are infinitely many primes~$r$ such that $\Omega(r^2-1)\leq 21$. This statement implies that $M_2\leq 21$, and, therefore, the following is a consequence of Theorem~\ref{t:counterexample1}.

\begin{cor} The $c$-dimensions of groups in $\mathscr{G}_1$ do not exceed $50$.
\end{cor}

It is clear that the number $p$ in the construction of the module $V_{1,r}$ can be chosen in such a way that $p$ would not divide $|X_{1,r}|$.
With this choice, one can apply the description of subgroups in $SL_2(r)\cong Sp_2(r)\cong A_{1,r}$, the character table of this group and results of \cite{Isaacs1} to calculate the character of the representation of $X_{1,r}$, which corresponds to the module $V_{1,r}$, and use it to calculate the precise value of $c$-dimension of $G_{1,r}$.

\sect{Proof of Theorem \ref{t:3Fitting}}

Throughout the rest of the paper $G$ denotes a finite group and $k$ denotes its $c$-dimension. The function $\lambda$ was defined in Section 1. We start by bounding the value of $\lambda$ on the layer $E(G)$ of $G$.

\begin{prop}\label{l:Quasisimple} There exists a universal constant $b$ such that ${\lambda(E(G))\leqslant b\cdot k}$ for every $G$.
\end{prop}

\begin{prf} Since the $c$-dimension of a subgroup does not exceed the $c$-dimension of the group, we may assume that $G=E(G)$. By Lemma \ref{l:CentralProduct}, we may suppose that $G$ is quasisimple. Denote by $\overline{\phantom{g}}$ the natural homomorphism from $G$ to $G/Z(G)$. Observe that if $\overline G$ is a sporadic group, a group of Lie type of bounded Lie rank, or an alternating group of bounded degree, then we can choose $b$ being large enough to make the proposition trivial in these cases. So we may assume that $\overline G$ is either an alternating group whose degree is large enough, or a classical group whose Lie rank is large enough.

Let us find a prime divisor $r$ of the order of $\overline G$ possessing the following two properties. First, it is coprime to the order of the center $Z(G)$. Second, the maximal length $l$ of a chain of nested centralizers of subsets of $r$-elements is bounded from below by a linear function of $\lambda(G)$. Assume that such $r$ is determined. It is clear that every subset of $r$-elements of $\overline G$ can be presented in the form $\overline{M}$ for some subset $M$ of $G$ also consisting of $r$-elements. Let $$M_1\subset M_2\subset\dots\subset M_l\subset G$$ be a chain of subsets of $r$-elements such that $$C_{\overline G}(\overline M_1)>C_{\overline G}(\overline M_2)>\dots>C_{\overline G}(\overline M_l).$$  By the properties of coprime action, $C_{\overline G}(\overline{M_i})=\overline{C_G(M_i)}$. Hence $$C_{G}(M_1)>C_{G}(M_2)>\dots>C_{G}(M_l).$$ Therefore $l$ is at most the $c$-dimension of $G$. Since it is bounded from below by a linear function of $\lambda(G)$, the proposition follows. Now it remains to determine such $r$ for every~$G$.

Consider the condition that $r$ does not divide the order of $Z(G)$. Every prime divisor of the order of $Z(G)$ is a prime divisor of the order of the Schur multiplier $M(\overline G)$ of $\overline G$. The orders of Schur multipliers of all finite simple groups can be found in \cite[Section 6.1]{GLS}. Due to these results, if $M(\overline G)$ is not a $\{2, 3\}$\nobreakdash-group, then $\overline G$ is isomorphic to $A_n(q)$ or $^2A_n(q)$ and any prime divisor of the order of $M(\overline G)$ divides $6(q^2-1)$. So $r$ should be chosen subject to these restrictions.

Consider the second condition on $r$, that is a linear bound on maximal length of a chain of centralizers of $r$-elements. If $G$ is an alternating group~$\operatorname{Alt}_n$, then $r$ can be chosen to be $5$. For $1\leqslant i\leqslant \frac{n}{5}$, we can choose nested sets $M_i$ consisting of $i$ disjoint $5$-cycles. Their centralizers form a strict chain whose length exceeds $\frac{n}{5}-1$. By definition of $\lambda$, this gives $\frac{\lambda(G)}{5}-1$ as a lower bound for $l$.

Let $\overline G$ be a classical group over a field of order $q$. According to \cite[Propositions 7--12]{Carter}, group $\overline G$ contains a subgroup isomorphic to a central product of quasisimple subgroups $H_1$, $\dots$, $H_s$ where every $H_i$ is a group of Lie type~$A_3$ or ${}^2A_3$ over a field of order $q$. Moreover, the number $s$ of these factors is at least $\frac{\lambda(G)-6}{4}$. By the Zsigmondy theorem, there exists an element $h_i$ of~$H_i$ of order coprime to $6(q^2-1)$. Hence $h_i$ is not central in $H_i$. Therefore, $$C_{\overline G}(h_1)>C_{\overline G}(h_1, h_2)>\dots>C_{\overline G}(h_1, h_2,\dots, h_s),$$ and we found a strict series of centralizers of subsets of $r$-elements whose length is at least $\frac{\lambda(G)-6}{4}$. The proposition is proved.

\end{prf}


\begin{prop}\label{l:CommonRank} There exists a universal constant $d$ such that $\lambda(G)\leqslant d\cdot k$ for every $G$.
\end{prop}

\begin{prf}Let $R$ be the soluble radical of $G$. If $P$ is a Sylow subgroup of $R$, then $G/R \cong N_G(P)/(R\cap N_G(P))$. So nonabelian composition factors of $N_G(P)$ and $G$ coincide. Furthermore, the $c$-dimension of $N_G(P)$ as a subgroup of~$G$ is at most $k$. Therefore, we may assume that $R$ is the Fitting subgroup of~$G$.

Put $\overline{G}=G/R$. Let $l$ be the number of nonabelian composition factors of the socle $\overline{L}$ of $\overline{G}$ (note that $\overline{L}$ is the direct product of nonabelian simple groups). The quotient $\overline{G}/\overline{L}$ is an extension of a soluble group by a subgroup of the symmetric group $\operatorname{Sym}_l$. Lemma~\ref{l:NumberOfCompFactors} implies that the number $l$ is less than $5k$. Hence $\lambda(\overline G/\overline L)<25k/4$ by Lemma \ref{CompFactorSn}. Thus it is sufficient to show that $\lambda(\overline L)\leqslant d'\cdot k$ for some $d'$. In particular, we may assume that $G$ coincides with the full inverse image of  $\overline{L}$ in $G$.

Define $\mathcal F$ to be the set of nonabelian composition factors of $C_G(R)$. For a prime $p$, denote by $\mathcal F_p$ the set of nonabelian composition factors of $G/C_G(O_{p'}(R))$ whose order is divisible by $p$ (here $O_{p'}(R)$ stands for the maximal $p'$-subgroup of $R$). By the Feit--Thompson theorem \cite{FT} and the Thompson--Glauberman theorem \cite[Chapter II, Corollary 7.3]{Glauberman}, the order of every finite nonabelian simple group is divisible by $2$ and is not coprime to $15$.  Since $R=O_{p'}(R)O_{q'}(R)$ for distinct primes $p$ and $q$, we have $C_G(O_{p'}(R))\cap C_G(O_{q'}(R))=C_G(R)$. Therefore, each nonabelian composition factor of $G$ is contained in $\mathcal  F\cup \mathcal  F_2\cup \mathcal  F_3\cup \mathcal  F_5$. So to prove the proposition, it is sufficient to to bound the sum of values of $\lambda$ in each of these four sets by a linear function in $k$.

Put $K=C_G(R)$. Since $\overline{K}=KR/R$ is normal in $\overline{G}$, it is a direct product of elements of $\mathcal F$. Take $\overline S$ in $\mathcal F$. Denote by $S$ the full inverse image of $\overline S$ in $K$. Then $S^{(\infty)}$, that is the least term in the derived series of $S$, is a perfect central extension of $\overline S$ which is normal in $K$, so it is a component of $K$. Therefore, all elements of $\mathcal F$ are composition factors of the layer $E(K)$. Since $E(K)$ is a subgroup of $E(G)$, the number $\lambda(K)$ is bounded by $b\cdot k$ for some constant~$b$ due to Proposition~\ref{l:Quasisimple}.

Put $K_p=G/C_G(O_{p'}(R))$ for $p\in\{2, 3, 5\}$. It is clear that $K_p$ is an extension of a direct product of some nonabelian composition factors of $G$ by a nilpotent $p'$-group. By the definition, $K_p$ acts faithfully on $O_{p'}(R)$, and so does its Sylow $p$-subgroup $P$. Group $P$ is the direct product of Sylow $p$-subgroups of elements of $\mathcal F_p$. By Lemma~\ref{l:ElemAbelian}, $P$ contains an elementary abelian $p$-group whose rank is bounded from below by $c\sum_{S\in\mathcal{F}_p} \lambda(S)$ for some positive constant~$c$. It follows from Lemma~\ref{LemmaKhukhro} that $c\sum_{S\in\mathcal{F}_p} \lambda(S)\leqslant k$.

Finally, $\lambda(G)\leqslant 3k/c+b\cdot k$ as required.
\end{prf}

Now we are ready to finish the proof of Theorem \ref{t:3Fitting}. Put $\overline G=G/F_3(G)$. Let $\overline F^*$ be the generalized Fitting subgroup of $\overline G$. The quotient $\overline G/\overline F^*$ is isomorphic to a subgroup of the group of outer automorphisms $\operatorname{Out}(\overline F^*)$ of~$F^*$. The latter is a subgroup of the direct product of $\operatorname{Out}(F(\overline G))$ and $\operatorname{Out}(E(\overline G))$. By Lemma~\ref{l:KhukhroTh}, the order of $F(\overline G)$ and, therefore, the order of $\operatorname{Out}(F(\overline G))$ is bounded in terms of $k$. Hence it is sufficient to prove that $\operatorname{Out}(E(\overline G))$ contains an abelian subgroup of $k$-bounded index.

If $E(\overline G)$ is a product of components $Q_1$, $Q_2$, $\dots$, $Q_s$, then $\operatorname{Out}(E(G))$ is a subgroup of a semidirect product of $\operatorname{Out}(Q_1)\times \operatorname{Out}(Q_2)\times\dots\times \operatorname{Out}(Q_s)$ and some subgroup of $\operatorname{Sym}_s$. According to Lemma~\ref{l:NumberOfCompFactors}, the number $s$ is bounded in terms of $k$, therefore, it is sufficient to prove that the direct product of the outer automorphism groups of the components possesses an abelian subgroup of $k$\nobreakdash-bounded index. Group $\operatorname{Out}(Q_i)$ is a subgroup of the outer automorphism group of the corresponding simple group $S_i=Q_i/Z(Q_i)$ \cite[Corollary 5.1.4]{GLS}. If $S_i$ is an alternating or sporadic group, then the order of $\operatorname{Out}(S_i)$ is at most~$4$. If $S_i$ is a group of Lie type, then the index of the cyclic subgroup of field automorphisms in $\operatorname{Out}(S_i)$ is bounded in terms of Lie rank of $S_i$ \cite[Theorems 2.5.1 and 2.5.12]{GLS}. Due to Proposition~\ref{l:CommonRank}, the Lie rank of every~$S_i$ is bounded in terms of~$k$. Define $A_i$ to be trivial, if $S_i$ is not a group of Lie type, and the group of field automorphisms of $S_i$ otherwise. The index of the direct product $A_1\times\dots\times A_s$ in $\operatorname{Out}(S_1)\times\dots\times \operatorname{Out}(S_s)$ is $k$-bounded by above arguments,. This completes the proof of the theorem.

\smallskip

Observe that since the group $\overline G/\operatorname{Soc}(\overline G)$ is a homomorphic image of the quotient group of $G/F_3(G)$ by its layer, the following is an immediate consequence of Theorem 2.

\begin{cor}\label{t:solubleRadical} Let $G$ be a finite group of $c$-dimension $k$ and $R$ be its soluble radical. Put $\overline G=G/R$. Then the group $\overline G/\operatorname{Soc}(\overline G)$ contains an abelian subgroup of $k$-bounded index.
\end{cor}

\begin{center}\textbf{Acknowledgments}\end{center}

The research was supported by RSF (project N 14-21-00065).

\end{document}